\documentclass[letterpaper, 10 pt, conference]{ieeeconf}  

\IEEEoverridecommandlockouts                              
\usepackage{booktabs}
\usepackage{makecell}
\usepackage{cite}
\usepackage{amsmath,amssymb,amsfonts}
\usepackage{algorithm2e}
\usepackage{graphicx}
\usepackage{textcomp}
\usepackage{xcolor}
\usepackage{comment}
\usepackage{ulem}
\usepackage{subcaption}
\usepackage{hyperref}
\usepackage{array}
\usepackage{multirow}
\usepackage{array}
\usepackage{booktabs}
\usepackage{amssymb}
\usepackage{pifont}
\newcommand{\cmark}{\ding{51}}
\newcommand{\xmark}{\ding{55}}

\newtheorem{theo}{Theorem}
\newtheorem{remarkEnv}{Remark}
\newenvironment{remark}[1][]{\begin{remarkEnv}}{\hfill$\blacklozenge$\end{remarkEnv}}
\newcommand{\Rho}{\mathrm{P}}

\newcommand{\data}{\mathcal{D}}

\newcommand{\obsf}{\hat{\rho}}

\newcommand{\truef}{\rho}

\newcommand{\disc}{\mathcal{X}}

\newcommand{\sensors}{\mathcal{M}}
\newcommand{\op}{\mathcal{G}}
\newcommand{\opapprox}{\op_{\theta}}
\newcommand{\system}{\mathcal{S}}

\newcommand{\systempidx}[3]{\hat{\system}^{#3}[#1,#2]}

\newcommand{\opobs}{\hat{\system}}
\newcommand{\opobsol}{\opobs^{ol}}
\newcommand{\opobsolr}{\opobs^{ol\text{-}r}}
\newcommand{\opobscl}{\opobs^{cl}}
\newcommand{\corr}{\mathcal{N}}
\newcommand{\corrapprox}{\corr_\psi}
\newcommand{\err}{\mathcal{E}}

\newcommand{\fspace}[2]{\mathcal{H}^{#2}_{#1}}
\newcommand{\rhospace}{\fspace{\Omega}{1}}
\newcommand{\opdomain}[1]{\fspace{\Omega}{#1}}

\newcommand{\kernelop}{\mathcal{K}}
\newcommand{\kernel}{\kappa}
\newcommand{\abs}[1]{\lvert#1\rvert}
\newcommand{\norm}[1]{\|#1\|}
\newcommand{\dtin}{\Delta t}
\newcommand{\state}{\Rho}
\newcommand{\stateest}{\hat{\Rho}}
\newcommand{\stateestt}[1]{\stateest_{#1}}
\newcommand{\stateestdelayed}[1]{\stateest_{#1 -\nout\Delta t}}

\newcommand{\stateestd}{\stateest^{\data}}
\newcommand{\stateesttd}[1]{\stateestd_{#1}}
\newcommand{\stateestdelayedd}[1]{\stateesttd{#1 -\nout\Delta t}}

\newcommand{\stateestu}{\stateest^{u}}
\newcommand{\stateesttu}[1]{\stateestu_{#1}}
\newcommand{\stateestdelayedu}[1]{\stateesttu{#1 -\nout\Delta t}}
\newcommand{\stateestdelayedincu}[1]{\stateesttu{#1 -\nout\Delta t + \Delta t}}

\newcommand{\etaest}{\stateest^{+}}
\newcommand{\etaestt}[1]{\etaest_{#1}}
\newcommand{\etaestdelayed}[1]{\etaestt{#1 -\nout\dtin}}

\newcommand{\etaestd}{\stateest^{\data +}}

\newcommand{\etaestu}{\stateest^{u +}}
\newcommand{\etaesttu}[1]{\etaestu_{#1}}
\newcommand{\etaestdelayedu}[1]{\etaesttu{#1 -\nout\dtin}}

\newcommand{\error}{\hat{e}^{+}}

\newcommand{\nout}{n_{d}}

\RestyleAlgo{ruled}
\SetKwComment{Comment}{// }{ }
\SetKwInOut{Input}{input}
\SetKwInOut{Output}{output}

\title{\LARGE \bf
Closed-Loop Neural Operator-Based Observer of Traffic Density}

\author{Alice Harting, Karl Henrik Johansson, and Matthieu Barreau
\thanks{This work is supported by the Wallenberg AI, Autonomous Systems and Software Program (WASP) funded by the Knut and Alice Wallenberg Foundation.}
\thanks{A. Harting, K.H. Johansson, and M. Barreau are with the Division of Decision and Control Systems, Digital Futures, KTH Royal Institute of Technology, SE-100 44 Stockholm, Sweden {\tt\small \{aharting,kallej,barreau\}@kth.se}}%
}
\begin{document}

\maketitle
\thispagestyle{empty}
\pagestyle{empty}

\begin{abstract}
We consider the problem of traffic density estimation with sparse measurements from stationary roadside sensors. Our approach uses Fourier neural operators to learn macroscopic traffic flow dynamics from high-fidelity data. During inference, the operator functions as an open-loop predictor of traffic evolution. To close the loop, we couple the open-loop operator with a correction operator that combines the predicted density with sparse measurements from the sensors. Simulations with the SUMO software indicate that, compared to open-loop observers, the proposed closed-loop observer exhibits classical closed-loop properties such as robustness to noise and ultimate boundedness of the error. This shows the advantages of combining learned physics with real-time corrections, and opens avenues for accurate, efficient, and interpretable data-driven observers.
\end{abstract}

\section{Introduction}
Freeway congestion is a major issue in metropolitan areas, leading to increased travel times and excessive fuel consumption \cite{Ferrara2018}. To address this, a core component of \textit{Intelligent Transportation Systems} is freeway traffic control, which operates at the macroscopic level (ramp management and variable speed limits) or at the vehicle level (intelligent vehicle-based control). Traffic control relies on measurements of state variables such as traffic density, which can be collected from stationary sensors (inductive loops) or mobile sensors (connected vehicles). However, these measurements are often sparse and noisy. To complete and denoise the state observation, estimation approaches combine data and prior knowledge of traffic dynamics \cite{trafficstateestimation}.

Model-based methods use physical models with calibrated parameters. Traffic state estimation then typically boils down to solving initial-boundary value problems (IBVPs) of partial differential equations (PDEs) with numerical methods such as the Godunov scheme \cite{lebacque1996godunov, barreau2020dynamic} or switching mode models \cite{smm}. To manage measurement noise and disturbance, it is common to use Kalman or particle filter-based approaches \cite{kf, mihaylova_particle_2004}. For example, \cite{bekiaris2016highway, bekiaris2017highway} reformulate the continuity equation into a linear parameter-varying (LPV) system, where real-time flow measurements are treated as parameters. While this class of methods typically requires little data and offers high interpretability, the accuracy is highly dependent on model assumptions and calibration \cite{trafficstateestimation}.

Data-driven methods, on the other hand, learn a system model directly from historical data, thereby relaxing explicit PDE assumptions and requirements of well-posedness \cite{trafficstateestimation}. Classical methods include ARIMA variants, Gaussian process (GP) models, and state-space Markov models. The rise of deep learning for nonlinear function approximation has enabled more expressive models of traffic flow prediction, with early examples such as the stacked autoencoder \cite{stackedAE}. In particular, Convolutional Neural Networks (CNNs) and Recurrent Neural Networks (RNNs) are often used as building blocks to capture spatiotemporal correlation. While deep learning approaches are expressive and require few assumptions, they are limited by a high data demand and computational complexity, and low interpretability. Thus, lightweight real-time prediction with explainable models remains an open challenge \cite{deeplearningtrafficsurvey}. 

In pursuit of data-efficient and interpretable models, hybrid models combine data-driven methods with prior physical knowledge. 
For example, \cite{BarreauMatthieu2021PLfI} reconstructs the traffic density state using physics-informed neural networks (PINNs) trained on offline data from sparse mobile sensors while regularizing with a continuous traffic PDE model. In \cite{yuan2022macroscopictrafficflowmodeling}, traffic models are used in physics-regularized GPs (PRGPs), which rely on latent force models and regularized kernel learning. However, these methods are restricted by explicit model assumptions on global traffic behavior. Furthermore, online estimation can be computationally demanding, since it requires retraining the network (PINNs) or performing matrix inversion with cubic complexity in terms of the size of the dataset (PRGPs).

We identify the need for a data-driven observer of traffic density that encodes a learned governing model in an interpretable manner, and offers lightweight real-time prediction with the integration of online measurements. For finite-dimensional nonlinear systems, there are learning-based solutions in the form of closed-loop data-driven observers where training is done offline and low-cost inference is performed online \cite{niazi2025kklobserversynthesisnonlinear}. The natural extension to PDE systems would be to consider Fourier neural operators (FNOs) \cite{li2021fourierneuraloperatorparametric}, which have emerged as a method to learn operators between infinite-dimensional function spaces from finite-dimensional data. 

\begin{table*}[t]
\centering
\renewcommand{\arraystretch}{1.3} 
\begin{tabular}{l ccc cc c}
\toprule
 &
\multicolumn{3}{c}{\textbf{Model approach}} &
\multicolumn{2}{c}{\textbf{Observation mechanism}} & \\
\cmidrule(lr){2-4} \cmidrule(lr){5-6} 
\textbf{Method} & \textbf{Model} & \textbf{Data} & \textbf{Learned closure} &
\textbf{Online inference} & \textbf{Closed-loop} & \textbf{Grid-free} \\
\midrule
Numerical integration of IBVP\cite{lebacque1996godunov, barreau2020dynamic, smm} & \cmark & \xmark & \xmark & \cmark & \xmark & \xmark \\
Filtered state-space model \cite{kf, mihaylova_particle_2004} & \cmark & \xmark & \xmark & \cmark & \cmark & \xmark\\
Filtered LPV-model with data-driven parameters \cite{bekiaris2016highway, bekiaris2017highway} & \cmark & \cmark & \xmark & \cmark & \cmark & \xmark\\
Stacked autoencoder \cite{stackedAE}; 
CNN + RNN\cite{deeplearningtrafficsurvey} & \xmark & \cmark & \cmark & \cmark & \xmark & \xmark \\
PINN \cite{BarreauMatthieu2021PLfI} & \cmark & \cmark & \cmark & \xmark & \cmark & \cmark \\ 
PRGP \cite{yuan2022macroscopictrafficflowmodeling} & \cmark & \cmark & \xmark & \cmark & \cmark & \cmark \\ 

KKL observer \cite{niazi2025kklobserversynthesisnonlinear}  & \cmark & \xmark & \xmark & \cmark & \cmark & \xmark \\
FNO \cite{fnomacro} & \cmark \& \xmark & \cmark \& \xmark & \cmark \& \xmark & \cmark & \xmark & \cmark \\
\textbf{FNO-based observer}& \cmark \& \xmark & \cmark \& \xmark & \cmark \& \xmark & \cmark & \cmark & \cmark \\
\bottomrule
\end{tabular}
\caption{Comparison of related methods for traffic state estimation.}
\label{tab:methods}
\end{table*}

Notably, \cite{fnomacro} successfully applied FNOs to learn a PDE solution operator for IBVPs and inverse problems of a first-order macroscopic model of traffic density flow. With this method, inference with new input data is significantly cheaper than with PINNs, since it is performed by a single forward pass through the operator network. However, the solver in this work is designed to generate either an offline reconstruction or open-loop prediction, without closed-loop integration of online data. 

Motivated by this, we introduce a data-driven, closed-loop neural observer of traffic density flow from sparse measurements. Our contributions are leveraging FNOs to learn a prediction operator of traffic density flow from high-fidelity data, and integrating Luenberger observer theory to get a robust density estimate using online measurements. We statistically evaluate the performance of the observer and compare it to open-loop variants. The relationship between the proposed observer and existing alternative methods is summarized in Table \ref{tab:methods}.

The paper is organized as follows. Section~\ref{sec:background} gives a brief introduction to traffic flow theory and neural operators, concluding with the problem statement. Section~\ref{sec:method} describes the methodology for designing three types of observers: open-loop, open-loop with reset, and closed-loop. Section~\ref{sec:resdisc} presents an extensive numerical evaluation of the proposed observers. Finally, Section~\ref{sec:conclusion} concludes the paper and outlines new research directions.



\textbf{Notation:} Let $\fspace{D}{d}\triangleq H^2_{\text{per}}(D; \mathbb{R}^d)$ denote the Sobolev space of order two containing periodic functions with domain $D$ and codomain dimension $d$. Further, we denote an operator with $\op$ or $\corr$. Finally, a GP is denoted $\mathcal{GP}$, and the posterior process associated with data $\data$ is denoted $\mathcal{GP}_\data$. We assume that the mean and covariance functions are chosen such that sample paths belong to $\fspace{D}{d}$ almost surely.

\section{Background and Problem Statement}\label{sec:background}
We begin by briefly introducing basic traffic modeling theory, and then present neural operators in this context. This is followed by a formal statement of the research problem.
\subsection{Modeling of Traffic Density Flow}\label{sec:modeling}
Traffic density is defined as the normalized number of vehicles per space unit at location $x$ and time $t$, $\rho(x,t)\in[0,1]$. Macroscopic models describe the dynamics of $\rho$ directly. While they are fast and simple, they typically rely on extensive assumptions. In contrast, microscopic models track each vehicle and offer greater precision, but are more computationally intensive and less convenient to use in freeway traffic control \cite{Ferrara2018}.
\subsubsection{Macroscopic scale}
A fundamental continuous first-order macroscopic model for traffic density flow is given by the \textit{Lighthill–Whitham–Richards} framework. It expresses conservation of vehicles through the one-dimensional continuity equation, where the flow $Q\in C^2([0,1])$ is uniquely determined by the local density \cite{lwr}:
\begin{equation}
\label{eq:lwr}
\begin{split}
    \frac{\partial \rho(x,t)}{\partial t} + \frac{\partial Q(\rho(x,t))}{\partial x} = 0&, \quad x \in \mathbb{R}, \quad t \geq 0,\\
    \rho(\cdot, 0)=\rho_0&, \quad x\in\mathbb{R}.
\end{split}
\end{equation} In general, classical solutions to the initial value problem of the hyperbolic PDE \eqref{eq:lwr} where $Q$ is smooth may develop discontinuities in finite time even when the initial data $\rho_0$ is smooth \cite{Dafermos}. The weak solutions, on the other hand, are non-unique but can be filtered for physical relevancy using an entropy condition associated with the conservation law \eqref{eq:lwr}. In fact, whenever $Q$ is smooth and $\rho_0$ is well-behaved, there exists a unique admissible (entropic) solution $\rho(\cdot, t)$ that is continuous in time and locally well-behaved \cite[Theorem 6.2.1, 6.2.2]{Dafermos}.
In this view, traffic density estimation is a well-defined problem. However, a key challenge for macroscopic models is to properly model the flow $Q$.

\subsubsection{Microscopic scale}To generate representations of higher fidelity without direct modeling of 
$Q$, traffic flow can instead be simulated using higher-order models at the vehicle level\cite{Ferrara2018}. These models integrate interactions between vehicles and road infrastructure to produce realistic traffic flows. Vehicle interactions are commonly represented by car-following models, in which the behavior relative to the preceding vehicle or platoon is specified. For example, \cite{krauss1998microscopic} proposed to let cars drive as fast as possible while maintaining a safe velocity to avoid collision. The model can be made increasingly complex by introducing parameters such as driver reaction time, estimation errors, and anticipation \cite{TREIBER200671}.    
In this way, microscopic models may generate high-fidelity traffic flows with few restrictive assumptions. However, a key challenge is the lack of a compact representation of the dynamics that effectively integrates with macroscopic-level data and control strategies. 
\subsection{Fourier Neural Operators (FNOs)}\label{sec:nos}

To address this challenge, neural operators provide a framework for learning the governing dynamics on the macroscopic level. Generally, neural operators learn a mapping between infinite-dimensional function spaces, that is, spaces containing functions with continuous domains, using finite-dimensional data. A special case is when the operator represents a forward solver for a hyperbolic PDE, such that the input function is the initial state and the output function is the solution at a final time $T$. 

We introduce neural operators in the context of traffic density flow arising from higher-order dynamics on a ring road of length $L$, denoted $\Omega\triangleq[0, L]$. The closed course induces realistic traffic flows, including stop-and-go waves formed by tiny fluctuations \cite{sugiyama2008traffic, tadaki2013phase}. Consider pairs of functions \[
(\Rho_0, \rho_T)\in\opdomain{n}\times\rhospace,
\] where $\Rho_0=[\rho(\cdot, 0), \rho(\cdot, -\dtin), \rho(\cdot, -2\dtin), ...]$ represents the discrete $n^{\text{th}}$-order initial state of step length $\dtin$, and $\rho_T= \rho(\cdot, T)$ represents the resulting traffic density at time $T$. Assume there exists an operator $\op$ such that
\begin{equation*}
    \op: 
\begin{array}[t]{rcl}
\opdomain{n}&\to&\rhospace, \\
\Rho_0&\mapsto&\rho_T,
\end{array}
\end{equation*} as can be motivated by the existence of a solution in the first-order case \eqref{eq:lwr}. Given observations of such pairs $\{(\Rho_0^i, \rho_T^i)\}_{i=1}^N$, where the initial states are generated according to some distribution $\Rho_0^i\sim\mu$, the goal is to learn an approximate solution operator $\opapprox\approx\op$ where $\opapprox$ is parametrized by a finite-dimensional vector $\theta$.

The approximate solution operator $\opapprox:\opdomain{n}\to\rhospace$ can be modeled by a neural operator \cite{li2021fourierneuraloperatorparametric}:
\begin{equation}\label{eq:fno}
\begin{split}
    v_1(x)=\mathcal{P}[\Rho_0](x),&\\
    v_{M}(x)=\mathcal{L}_{M}\circ ...\circ\mathcal{L}_1[v_1](x),&\\
    \rho_T(x)=\mathcal{Q}[v_{M}](x),\quad\forall x\in \Omega,&
\end{split}
\end{equation}
where $\mathcal{P}$ is defined by a local transformation $P:\mathbb{R}^{n}~\times~\Omega\to\mathbb{R}^{d_v}$, $\mathcal{P}[\Rho_0](x)\triangleq P(\Rho_0(x), x)$ acting on function inputs and outputs directly to lift the input state to higher-dimensional space of size $v$. Similarly, $\mathcal{Q}$ is defined by a local transformation $Q:\mathbb{R}^{d_v}\times\Omega\to\mathbb{R}$, $\mathcal{Q}[v](x)\triangleq Q(v(x),x)$ projecting back to the physical space. Both are typically shallow neural networks. The intermediate \textit{spectral convolution layers} $\mathcal{L}_l:\mathcal{V}\to\mathcal{V}$, where $\mathcal{V}=\mathcal{V}(\Omega;\mathbb{R}^{d_v})$, are structured as 
\begin{align*}
    \mathcal{L}_l[v](x)=\sigma\big(W_l[v](x)+\kernelop_{\phi_l}[v](x)\big),\quad\forall x\in \Omega,
\end{align*}
where $\sigma$ is a local nonlinear activation function, and $W_l$ is a local pointwise affine transformation. The key component for the infinite-dimensional setting is the non-local bounded linear operator $\kernelop_{\phi_l}\in\mathcal{L}(\mathcal{V}, \mathcal{V})$, which acts in function space.

A central aspect of FNOs \cite{li2021fourierneuraloperatorparametric} is that $\kernelop_{\phi_l}$ is defined by the convolution 
\begin{equation}
\label{eq:conv}
    \kernelop_{\phi_l}[v]\triangleq\kappa_{\phi_l} * v,
\end{equation}
where $\kernel_{\phi_l}:\Omega \to\mathbb{R}^
{d_v\times d_v}$ is a periodic, integrable function which has a Fourier series representation\footnote{This aligns with our problem setting, as the ring road formulation entails that $\rho(\cdot, t)$ is $L$-periodic on $\Omega$.}. By the convolution theorem, \eqref{eq:conv} admits the frequency domain form
\begin{align}
    \kernelop_{\phi_l}[v](x) = \mathcal{F}^{-1}\big(\mathcal{F}(\kappa_{\phi_l})(k)\cdot \mathcal{F}(v)(k)\big)(x),
    \label{eq:kernellayer}
\end{align} where $\mathcal{F}$ computes the Fourier coefficients and $k\in\mathbb{Z}$. Parameterizing the convolution kernel directly in the frequency domain, $R_{\phi_l}\triangleq\mathcal{F}(\kappa_{\phi_l})$, and truncating the Fourier series to $k_{max}$ frequency modes, $R_{\phi_l}\in\mathbb{C}^{k_{max}\times d_v\times d_v}$ and $\mathcal{F}(v)\in\mathbb{C}^{k_{max}\times d_v}$, yields a finite-dimensional operator parametrization. Despite this truncation, the operator \eqref{eq:conv} remains defined on functions over the continuous domain $\Omega$, since $\mathcal{F}^{-1}$ reconstructs from the basis functions $e^{2\pi i\langle \cdot,k\rangle}$ for any $x\in \Omega$. In practice, input–output data consist of discretized rather than analytic functions, and the Fourier transform $\mathcal{F}$ can be computed using the Fast Fourier Transform $\hat{\mathcal{F}}$. Crucially, truncating in the frequency domain -- rather than the physical domain -- allows for a parametrization that can be shared across resolutions. In other words, the resulting operator \eqref{eq:fno}--\eqref{eq:kernellayer} is discretization convergent, and thus consistent in the infinite-dimensional setting.



Similar to neural networks, FNOs are capable of approximating a large class of operators between infinite-dimensional spaces up to arbitrary accuracy. Below, we present the universal approximation theorem from \cite{kovachki2021universalapproximationerrorbounds} as applied to the setting of this paper.
\begin{theo} 
    \noindent Let $\op:\opdomain{n} \to \opdomain{1}$ be a continuous operator, and $K\subset \opdomain{n}$ a compact subset. Let $\opapprox:\opdomain{n}\to \opdomain{1}$ be defined as in \eqref{eq:fno}--\eqref{eq:kernellayer} with the coefficients $\theta\triangleq\{\theta_P, \theta_Q, \theta_{W_l}, R_{\phi_l}\}$. Then, $\forall\epsilon>0$ $\exists  k_{max}\in\mathbb{N}$ and $\theta\in\Theta$ such that $\opapprox$, continuous as an operator, satisfies
\begin{equation}
\label{eq:universalapprox}
    \sup_{a\in K} \int\limits_\Omega\norm{\op[a](x)-\opapprox[a](x)}^2dx\leq \epsilon.
\end{equation}
\end{theo}It is further suggested by \cite{kovachki2021universalapproximationerrorbounds} that, while the approximation error scales super-exponentially with the size of the FNO in the general case, it typically scales polynomially for PDEs of the form \eqref{eq:lwr} under suitable assumptions. The expression \eqref{eq:universalapprox} serves as the basis for a training loss function, and may be extended to include PINN-type regularization from a predefined PDE model \cite{li2024physics}.
\subsection{Problem Statement}\label{sec:ps} 
Following this discussion, we model traffic density flow as an $n^{\text{th}}$-order discrete-time dynamical system 
\begin{align}\label{eq:sys}
    \system[\bar{\Rho}_0]:
    \begin{cases}
        \truef(\cdot, t+\dtin) = \op[\state_t],\quad t > 0, \\
        \state_t=[\rho(\cdot, t), \rho(\cdot, t-\dtin), \dots],\\
        \state_{0}=\bar{\Rho}_0,
    \end{cases}
\end{align} 
where $\truef(\cdot, t) \in \rhospace$ is the true density,  $\state_t\in\opdomain{n}$ is the state at time $t$, $\bar{\Rho}_0$ is the initial state, and $\op$ represents the unknown solution operator of the dynamics.

\begin{figure}
    \centering
    \includegraphics[width=\linewidth]{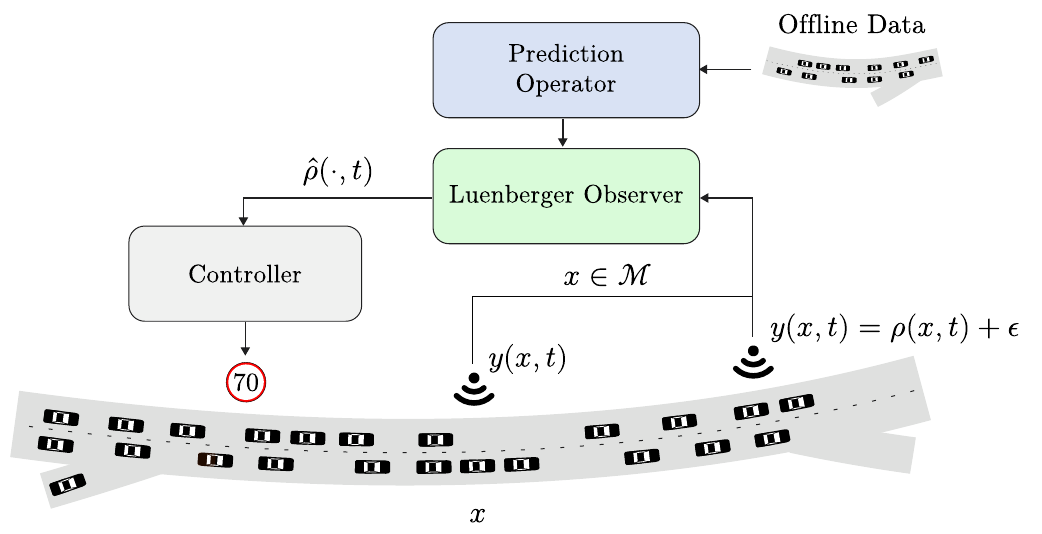}
    \caption{Problem setting. Our goal is to design an observer using offline data to learn the dynamics, and online measurements from stationary sensors at $x\in\mathcal{M}$ to calibrate predictions. The observer may then be used for macroscopic-level traffic control.}
    \label{fig:vehicle-road}
\end{figure}

The objective of this paper is to develop an algorithm in the form of a dynamical system that observes \eqref{eq:sys} in closed-loop by generating real-time predictions $\hat{\rho}$ of the system state $\rho$ every $\Delta t$ time unit such that $\| \rho - \hat{\rho}\|$ is minimized. To achieve this goal, we assume access to online measurements $y(x,t)=\rho(x,t)+\epsilon$ with i.i.d. noise $\epsilon$ at stationary, sparse sensor locations $x\in\sensors\subset\Omega$, collectively denoted $\mathbf{y}_{\sensors}(t)$. We also assume access to offline observations $\{\system[\Rho_0^i]\}_{i\in\mathcal{I}}$, where $\Rho_0^i\sim\mu$ for some distribution $\mu$. The problem setting is illustrated in Fig. \ref{fig:vehicle-road}.

\begin{remark}
    The problem is inherently ill-posed due to the incomplete observation of the initial state. Thus, our goal is to learn the most probable solution that remains physically consistent with the dynamics in the offline observations.
\end{remark}

\section{Methodology}\label{sec:method}

To build a closed-loop observer of \eqref{eq:sys}, the first step is to obtain a representation of the open-loop prediction, $\op$. Thus, our approach starts by learning an approximate solution operator $\opapprox \approx \op$ from offline data. We then introduce two methods of iteratively applying $\op$ to observe \eqref{eq:sys} in an open-loop manner, $\opobsol$ and $\opobsolr$. To close the loop, we finally design an algorithm, $\opobscl$, that integrates sparse online measurements with recursive predictions using $\op$. The research problem is then addressed by using $\opobscl$ together with $\opapprox$, resulting in a data-driven, closed-loop observer.

\subsection{Solution Operator: Model Identification}\label{sec:solop}
Given the offline scenarios $\{\system[\Rho_0^i]\}_{i\in\mathcal{I}}$, the goal is to learn a representation of the solution operator
\begin{equation*}
    \op_k: 
\begin{array}[t]{rcl}
\opdomain{n}&\to&\rhospace, \\
\Rho_0&\mapsto&\rho_{k\dtin},
\end{array}
\end{equation*}
which maps the input state $\Rho_0$ to the solution function at $k$ time steps later, $\rho_{k\dtin}$. While both $\Rho_0$ and $\rho_{k\dtin}$ are infinite-dimensional functions, we assume that the available data is represented in finite-dimensional form, obtained by evaluating them on a uniformly discretized grid $\disc\subset\Omega$. This yields $\Rho_0\rvert_\disc\in\mathbb{R}^{n\times \abs{\disc}}$ and  $\rho_{k\dtin}\rvert_\disc\in\mathbb{R}^{\abs{\disc}}$. To bridge the gap between finite-dimensional data and function space learning, the approximate operator $\opapprox$ is modeled by an FNO \eqref{eq:fno}--\eqref{eq:kernellayer}. 

\begin{figure}[t!]
    \centering
    \includegraphics[width=\linewidth]{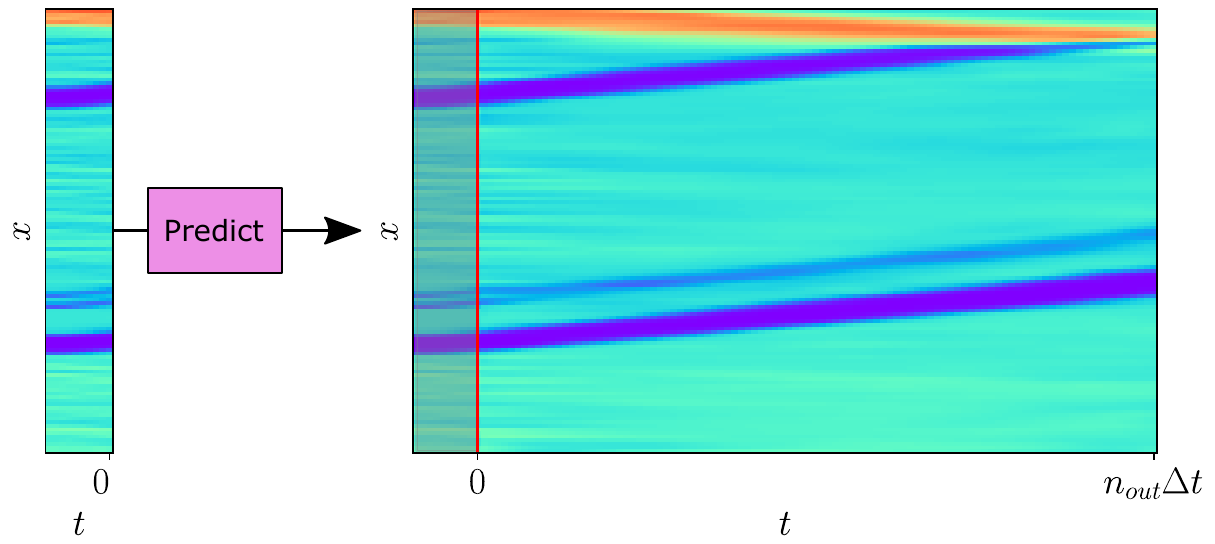}
    \caption{Open-loop prediction with solution operator $\opapprox$. It predicts $\obsf_{\dtin},\obsf_{2\dtin},...,\obsf_{n_{out}\dtin}$ given an initial state $\state_0=[\rho_{0}, \rho_{-\dtin},\dots,\rho_{-(n-1)\dtin}]$. Here, $n=10$ and $n_{out}=100$.}
    \label{fig:architechture-ol}
\end{figure}

With a target offset $n_{out}$ for the observer predictions, $\opapprox$ is trained to predict the density over the entire horizon,
\begin{equation}\label{eq:opdef}
\opapprox:
\begin{array}[t]{rcl}
\opdomain{n} & \to & \opdomain{n_{out}}, \\
\Rho_0 & \mapsto & [\obsf_{\dtin}, \obsf_{2\dtin}, \dots, \obsf_{n_{out}\dtin}].
\end{array}    
\end{equation}
Including the horizon up to the target estimate $\obsf_{n_{out}\dtin}$ was empirically found to improve performance, and it reflects the setup in \cite{gopakumar2024uncertaintyquantificationsurrogatemodels}. Thus, the corresponding data-driven learning objective\footnote{While we focus on a data-driven objective in this paper, we highlight that it is possible to regularize using lower-order models such as \eqref{eq:lwr} with a predefined or jointly learned flow function $Q$ following \cite{li2024physics, fnomacro}.} is given by
\begin{equation}\label{eq:solop}
\min_{\theta}\mathbb{E}_{a\sim\mu}\bigg[\sum\limits_{k=1}^{n_{out}}\int\limits_\Omega \norm{\op_{k}[a](x)-\opapprox^{(k)}[a](x)}^2dx\bigg].
\end{equation}Fig. \ref{fig:architechture-ol} illustrates an example forward pass of a learned operator $\opapprox$ for $n_{out}=100$. We refer to this as open-loop prediction in the observer framework.
\subsection{Open-Loop Observers}\label{sec:ol}
Given a representation of $\op$, the goal is now to develop the observer algorithm. We begin by considering standard autoregression, wherein past estimates serve as input to succeeding predictions through the \textit{predicted state estimate} 
\begin{equation*}
    \stateestt{t} = [\obsf(\cdot, t), \obsf(\cdot, t-\dtin), \dots]\in\opdomain{n}.
\end{equation*} Thus, triggering an autoregressive rollout only requires an estimate of the initial state $\stateestt{0}$. However, this is only partially observed in our problem, since measurements are restricted to the sparse sensor locations $\sensors\subsetneq\Omega$. To complete the initial state estimate, we interpolate between the sparse measurements $\mathbf{y}_{\sensors}(t)$ with GP regression,
\begin{equation*}
    \stateestt{0}\sim\mathcal{GP}_{[\mathbf{y}_{\sensors}(0), \mathbf{y}_{\sensors}(-\dtin),\dots]} \in\opdomain{n}.
\end{equation*}Thereafter, the estimates $\obsf(\cdot,t)$ are generated by shifting the autoregressive prediction window in increments of $\dtin$. 

\begin{figure}[ht!]
    \centering
    \includegraphics[scale=0.45]{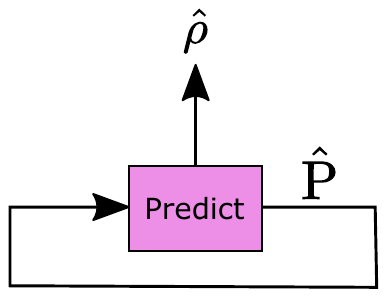}
    \caption{Open-loop observer $\opobsol$}
    \label{fig:open-loop}
\end{figure}

We found that using multi-step ahead prediction $\op_{\nout + 1}$ with a delayed state estimate $\stateestdelayed{t}$ offset by $\nout>n-1$ improved stability. This creates the \textit{open-loop observer} $\opobsol$:
\begin{align}\label{eq:obsol}
    \systempidx{\op}{\mathbf{y}_\sensors}{ol}:
    \left\{
    \begin{array}{l}
        \obsf(\cdot, t+\dtin) = \op_{\nout + 1}[\stateestdelayed{t}], \quad \quad \quad \\
        \hfill t = 0,\dtin,\dots\\
        \stateestdelayed{t}=[\obsf(\cdot, t-\nout\dtin), \dots],\\
        \stateestdelayed{ }\sim\mathcal{GP}_{[\mathbf{y}_{\sensors}(-\nout\dtin),\dots]}.
    \end{array}
    \right.
\end{align} The gap between the initial state and the first estimate, $\obsf(\cdot, t),$~$t\in[-\nout\dtin + \dtin,\dots, 0]$, is filled by interpolating the corresponding measurements with GP regression. We illustrate the observer in Fig. \ref{fig:open-loop}.


While the purely autoregressive approach benefits from consistent predictions, it is sensitive to disturbance in the initial state. As reported in \cite{mccabe2023towards}, autoregressive rollouts of neural operators suffer from unbounded out-of-distribution (OOD) error growth for long-term forecasts. Since we have access to online measurements, we compare $\opobsol$ with another open-loop observer that resets the input state for every new prediction based on recent measurements. Specifically, the \textit{data-based state estimate} $\stateesttd{t}$ is obtained by interpolating the corresponding measurements $[\mathbf{y}_{\sensors}(t), \mathbf{y}_{\sensors}(t-\dtin),\dots]$ using GP regression,
\begin{equation}
    \stateesttd{t}\sim\mathcal{GP}_{[\mathbf{y}_{\sensors}(t), \mathbf{y}_{\sensors}(t-\dtin),\dots]} \in\opdomain{n}.
\end{equation}

\begin{figure}[ht!]
    \centering
    \includegraphics[scale=0.45]{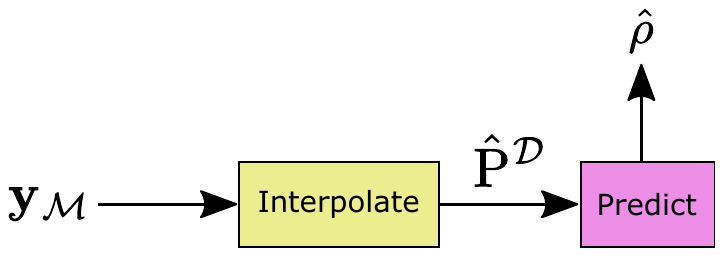}
    \caption{Open-loop observer with reset $\opobsolr$}
    \label{fig:open-loop-reset}
\end{figure}

For comparability with \eqref{eq:obsol}, we consider a delayed state estimate $\stateestdelayedd{t}$. The resulting dynamical system is denoted \textit{open-loop observer with reset} $\opobsolr$, and defined as
\begin{align}
    \label{eq:obsolr}
        \systempidx{\op}{\mathbf{y}_\sensors}{ol\text{-}r}:
    \left\{
    \begin{array}{l}
        \obsf(\cdot, t+\dtin) = \op_{\nout + 1}[\stateestdelayedd{t}], \quad \quad \quad \\
        \hfill t = 0, \dtin,\dots\\
        \stateestdelayedd{t}\sim\mathcal{GP}_{[\mathbf{y}_{\sensors}(t - \nout\dtin),\dots]},\\
        \stateestdelayedd{}\sim\mathcal{GP}_{[\mathbf{y}_{\sensors}(-\nout\dtin),\dots]}.
    \end{array}
    \right.
\end{align} The observer is illustrated in Fig. \ref{fig:open-loop-reset}.

A limitation of $\opobsolr$ is that each prediction relies on a few data points. By the property of the disconnected predictions, the observer does not leverage the full historical dataset.
\subsection{Closed-Loop Observer}\label{sec:cl}
\begin{figure}[ht!]
    \centering
    \includegraphics[scale=0.45]{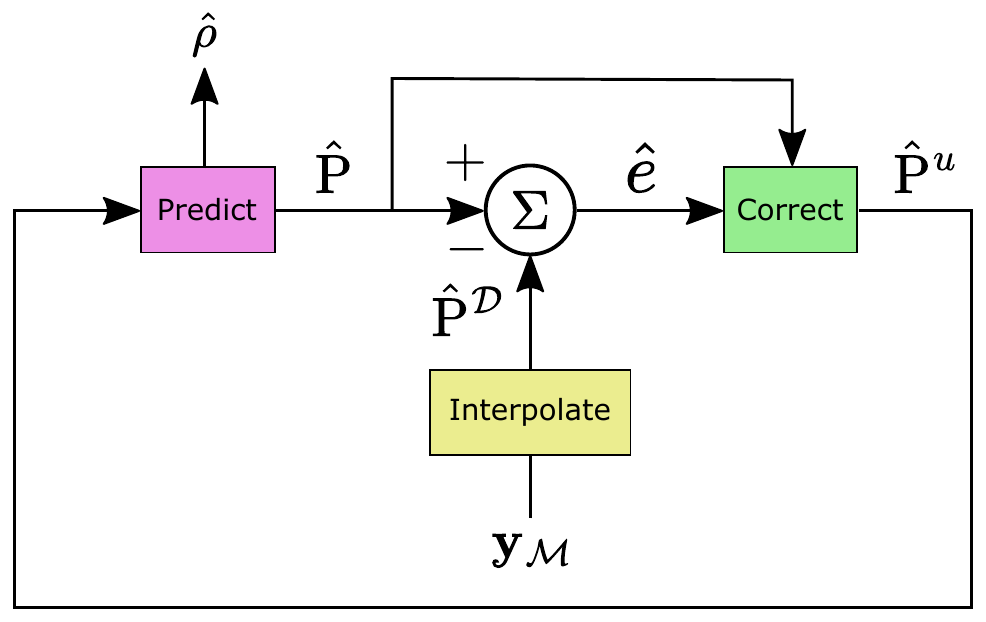}
    \caption{Closed-loop observer $\opobscl$}
    \label{fig:closed-loop}
\end{figure}

Combining the consistent predictions of $\opobsol$ with the online updates of $\opobsolr$ leads to our main contribution: a \textit{closed-loop observer} $\opobscl$ that integrates an autoregressive rollout of $\op$ with online measurements using a correction operator $\corr$. Analogous to a Kalman filter, $\corr$ adjusts the predicted state estimate using recent measurements before it is passed as input to the next prediction. Measurements are introduced via a Luenberger-type error operator $\err$, which computes the error estimate $\hat{e}$ as the difference between the predicted state estimate $\stateest$ and the data-based estimate $\stateestd$. The correction operator $\corr$ then processes $\stateest$ and $\hat{e}$ to generate an \textit{updated state estimate} $\stateestu\in \opdomain{n}$, which is passed as input to the next prediction. We summarize $\opobscl$ as
\renewcommand{\arraystretch}{1.25} 
\begin{align}
\label{eq:obscl}
    \hspace*{-0.2cm}
    \systempidx{\op, \corr}{\mathbf{y}_\sensors}{cl}:
    \left\{
    \begin{array}{l}
        \obsf(\cdot, t+\dtin) = \op_{\nout + 1}
        [\stateestdelayedu{t}], \quad \quad \quad \\
        \hfill t = 0, \dtin,\dots\\
        \stateestdelayedu{t}=\corr\big[\stateest, \hat{e}\big],\\
        \hat{e} = \err[\stateest, \stateestd],\\
        \stateestdelayed{t} = [\obsf(\cdot, t-n_d\dtin), \dots],\\
        \stateestdelayedd{t}\sim\mathcal{GP}_{[\mathbf{y}_{\sensors}(t - n_d\dtin),\dots]},\\ 
        \stateestdelayedu{ }\sim\mathcal{GP}_{[\mathbf{y}_{\sensors}(-\nout\dtin),\dots]},
    \end{array}
    \right.
\end{align}and illustrate the system in Fig. \ref{fig:closed-loop}.

\begin{figure}[ht!]
    \centering
    \includegraphics[width=\linewidth]{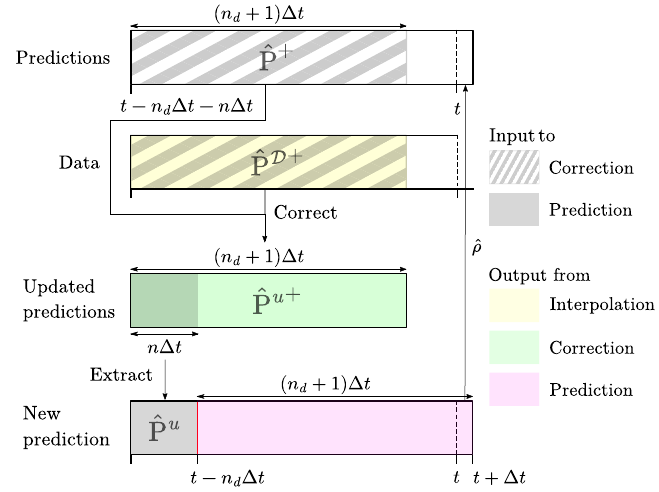}
    \caption{Implementation of $\opobscl$ with extended state estimates}
    \label{fig:extended_state_estimate}
\end{figure}

In implementing \eqref{eq:obscl}, corrections are applied to a moving window of historical predictions,
\begin{equation}\label{eq:etadef}
    \begin{array}{l}
            \displaystyle \etaestdelayed{t} = [\obsf(\cdot, t - (n-1)\dtin),\dots,\quad\quad\quad\quad\\
            \hfill\obsf(\cdot, t-\nout\dtin),\dots,\quad\quad\quad\quad\\
            \hfill\obsf(\cdot, t - \nout\dtin - (n-1)\dtin)],\\
    \end{array}
\end{equation} assuming $\nout>n-1$ as in the previous section. Thus, the \textit{extended} state estimate $\etaest$ contains $\stateest$ and extends forward in time such that the window length matches the prediction horizon \eqref{eq:opdef} where  $n_{out}=\nout+~1$ similar to \eqref{eq:obsol} and \eqref{eq:obsolr}. After applying the correction to this window, the updated state estimate $\stateestu$ is taken as the subset $\etaestu_{-n:-1}$. The procedure is illustrated in Fig. \ref{fig:extended_state_estimate}.

With $\etaest$ and the corresponding extended data-based estimate $\etaestd$, the Luenberger error operator $\err$ is defined as
\begin{equation*}
\err[\etaest, \etaestd]\triangleq \etaest - \etaestd\in \opdomain{\nout + 1}.  
\end{equation*}The correction operator $\corr$ processes $\etaest$ and the error $\error = \err[\etaest, \etaestd]$ to generate the updated state estimate $\etaestu$,
\begin{equation}
\label{eq:corrdef}
    \corr: 
    \begin{array}[t]{ccl}
    \opdomain{n_d+1} \times \opdomain{n_d+1} & \to & \opdomain{n_d+1}, \\
    \etaest, \error & \mapsto & \etaestu.
    \end{array}
\end{equation} The correction update window $\etaestdelayed{t}$ is shifted in increments of $\dtin$ in parallel with the prediction input window $\stateestdelayed{t}$, cf.~\eqref{eq:obsol}. Since $\stateestdelayedu{t}\subset\etaestdelayedu{t}$, the state estimates are updated before being passed as input to the next prediction.
The implementation details are outlined in Algorithm \ref{alg:clobs}.
\renewcommand{\arraystretch}{1} 
\begin{algorithm}[hbt!]
\caption{Closed-Loop Observer $\opobscl$}\label{alg:clobs}
\Input{$\op, \corr, \mathbf{y}_{\sensors}$}
\Output{$\obsf$}
$t=0$\;
$\obsf(\cdot, t)\sim\mathcal{GP}_{\mathbf{y}_{\sensors}(t)}, \quad -\nout\dtin<t\leq 0,\quad x\in\disc$\;
$\stateestdelayedu{}\sim\mathcal{GP}_{[\mathbf{y}_{\sensors}(-\nout\dtin),\dots]},\quad x\in\disc$\;
\Repeat{end}{
\Comment{Predict}
$\obsf(\cdot, t+\dtin) \gets \op_{\nout+1}[\stateestdelayedu{t}],\quad x\in\disc$\;
$\etaest \gets [\obsf(\cdot, t-(n-1)\dtin + \dtin), \dots]$\;
\Comment{Collect measurements}
$\etaestd \sim\mathcal{GP}_{[\mathbf{y}_{\sensors}(t - (n-1)\dtin+\dtin),\dots]}$\;
\Comment{Correct}
$\error=\err[\etaest, \etaestd]$\;
$\etaestu \gets \corr\big[\etaest, \error\big]$\;
\Comment{Extract next input}
$\stateestdelayedincu{t}\gets \etaestu_{-n:-1}$\;
$t\gets t+\dtin$\;}
\Return{$\obsf$}
\end{algorithm}


The correction operator $\corr$ is approximated using an FNO defined on a two-dimensional domain:
\begin{equation}
    \corrapprox:\fspace{\Omega\times\Gamma}{1}\times \fspace{\Omega\times\Gamma}{1}\to \fspace{\Omega\times\Gamma}{1}, 
\end{equation}
 where $\Gamma=[\dtin, (\nout + 1)\dtin]$ and the processed functions are translated in time to take values in this domain.
Note that $\Gamma$ reflects the choice of window size in \eqref{eq:etadef}. This setting enables an extension of the learning objective for open-loop prediction \eqref{eq:solop}, such that $\corrapprox$ is obtained as the minimizer of 
\begin{equation}\label{eq:corr}
\hspace*{-0.4cm}
\begin{array}{l}
    \displaystyle \min_{\psi}\mathbb{E}_{a\sim\mu}\bigg[\sum\limits_{k=1}^{\nout+1}\int\limits_{\Omega} \norm{J_k[a](\psi; x)}^2dx\bigg], \quad \quad \quad \quad \\ \\
    \displaystyle J_{k}[a](\psi; \cdot) \triangleq \op_{k}[a](\cdot) - \\\hfill\corrapprox\big[\opapprox[a], \err[\opapprox[a], \rho_{\data}[a]]\big](\cdot, k\dtin),\\
\end{array}
\end{equation} where $\rho_{\data}[a]$ is a data-based estimate of the true density $\mathcal{G}_{1:\nout+1}[a]$ using samples at the sensor locations, $\{\mathcal{G}_{1:\nout+1}[a](\cdot)\}_{x\in \sensors}$. By comparing the objective function of $\corrapprox$ in \eqref{eq:corr} with that of the open-loop predictor $\opapprox$ in \eqref{eq:solop}, we observe that $\corrapprox$ is trained to refine the prediction of $\opapprox$. This refinement is obtained by taking in measurements in a Luenberger manner through $\err$. 
\section{Numerical Simulations and Discussion}\label{sec:resdisc}
Numerical simulations are conducted using the open-source microscopic traffic simulator SUMO \cite{sumo}. This framework integrates road network data, infrastructure like traffic lights, and stochastic demand to produce realistic vehicle flows. Following the methodological section, we focus on ring road dynamics \cite{MatthSUMO}: vehicles are gradually introduced to the road during an initialization phase to reach the initial state, after which the higher-order microscopic dynamics are run in a loop. To convert the discrete distribution of vehicles to a continuous density function, the grid of vehicle-populated cells is convolved with a Gaussian filter \cite{sumoconv}.

The remainder of this section is organized as follows: we first outline implementation details, then demonstrate the observers on a single test case, evaluate their performance in noiseless, noisy, and out-of-distribution (OOD) scenarios, and finally analyze input-to-state stability.
\subsection{Implementation Details}\label{sec:impl}
The SUMO simulations are generated on a circular road of length $L=6.2$ km, discretized into 123 cells $\disc$ with six equidistant sensors $\sensors$. The simulations are run for a total time of $T=40$ minutes divided into time steps of length $\Delta t=1$ s, following \cite{MatthSUMO}. The training dataset contains 20 noiseless, independent simulations per average $\rho\in[0.1, 0.2, ..., 0.8]$ split into a total of 3360 non-overlapping input-output pairs $\big(\Rho_0^j, [\rho^j_{k\dtin}]_{k=1}^{\nout + 1}\big)$. We found that $n=10$ and $\nout+1=100$ works well in this setting. 

The test dataset consists of 10 simulations per average $\rho\in[0.3, ..., 0.8]$, focusing on the settings where congestion occurs. In addition, we test the observers on a noisy test set and a noiseless OOD test set. The noisy data is produced by adding i.i.d. Gaussian noise with standard deviation $\sigma=0.1$ to the measurements. The OOD data is generated by varying parameters in SUMO to induce more traffic jams.

The approximate solution operator $\opapprox$ and correction operator $\corrapprox$ share a similar FNO structure: the lifting function is a linear layer $P$ of output dimension 16, and the projection function $Q$ is a fully-connected neural network $Q$ with one hidden layer of size 128. In between, $\opapprox$ has four spectral convolution layers $\mathcal{L}_{l}$ where $W_l$ are one-dimensional convolution layers of sizes $[24, 24, 32, 32]$ and the Fourier expansions are truncated to $k_{\max}$ modes $[15, 12, 9, 9]$. The correction operator $\corrapprox$ has two spectral convolution layers $\mathcal{L}_{l}$ where $W_l$ are two-dimensional convolution layers of sizes $[24, 32]$ and the Fourier expansions are truncated to $k_{\max}$ modes $[15, 9]$. For both operators, we use GELU activation functions and warp the output with a sigmoid to impose $\obsf\in[0, 1]$. GP regression is implemented using a prior with a zero mean function and a squared exponential kernel with a length scale of 1, and we sample once during training while taking the posterior mean function during inference.

We use the Adam optimizer and train for 500 epochs. The training and testing were run on a laptop with an Intel(R) Core(TM) i7-1370p CPU @ 1.90 GHz and 14 processing cores. The code and experiments are available on GitHub for reproduction\footnote{\href{https://github.com/aharting/closed-loop-neural-operator-based-observer-of-traffic-density}{https://github.com/aharting/closed-loop-neural-operator-based-observer-of-traffic-density}}. 
\subsection{Numerical Simulations and Discussion}\label{sec:numsim}
\begin{figure}
    \centering
    \includegraphics[width=\linewidth]{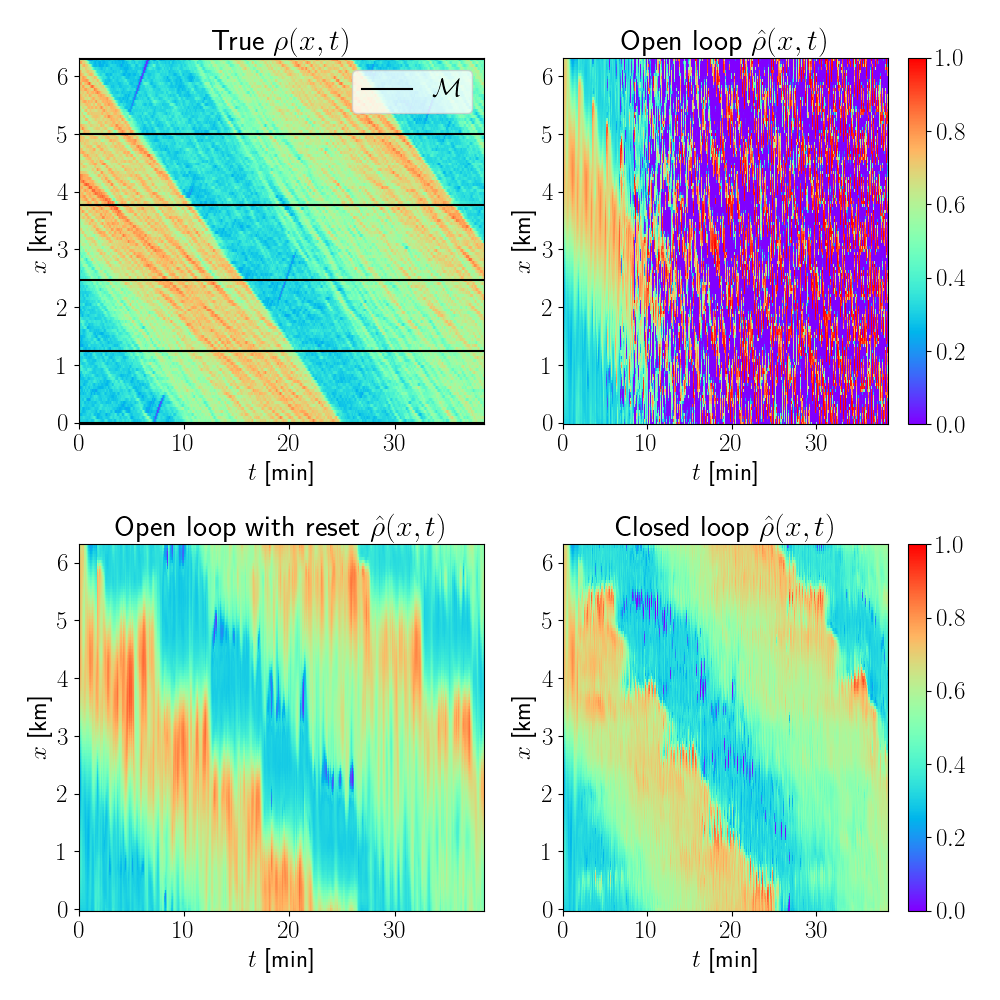}
    \caption{Traffic density estimation with observer variants. The true density is shown in the top left figure, where measurement locations are indicated in black. The density estimates are shown for the open-loop observer (top right), the open-loop observer with reset (bottom left), and the closed-loop observer (bottom right). Note that the solution periodicity is a result of the ring road implementation.}
    \label{fig:illustration_full_unravel}
    \vspace*{-0.5cm}
\end{figure}

A comparison between the open-loop observer $\opobsol$, open-loop observer with reset $\opobsolr$, and closed-loop observer $\opobscl$ applied to a single test example is shown in Fig. \ref{fig:illustration_full_unravel}. We observe that $\opobsol$ quickly diverges, while $\opobsolr$ produces disconnected predictions, with the vertical dispersion of the measurements resulting from the vertical GP interpolation between observation points. In contrast, $\opobscl$ maintains stability while producing consistent predictions over the spatiotemporal domain. This highlights the benefits of $\opobscl$: the continuous calibration of the autoregressive inputs both improves the prediction estimates and counteracts instability. Hence, $\opobscl$ unites the consistent predictions of $\opobsol$ with the online updates of $\opobsolr$, achieving an improved observer. We continue by exploring the performance across the test set in terms of accuracy, robustness, and input-to-state stability.

\begin{figure}[ht!]
    \centering
    \includegraphics[scale=0.28]{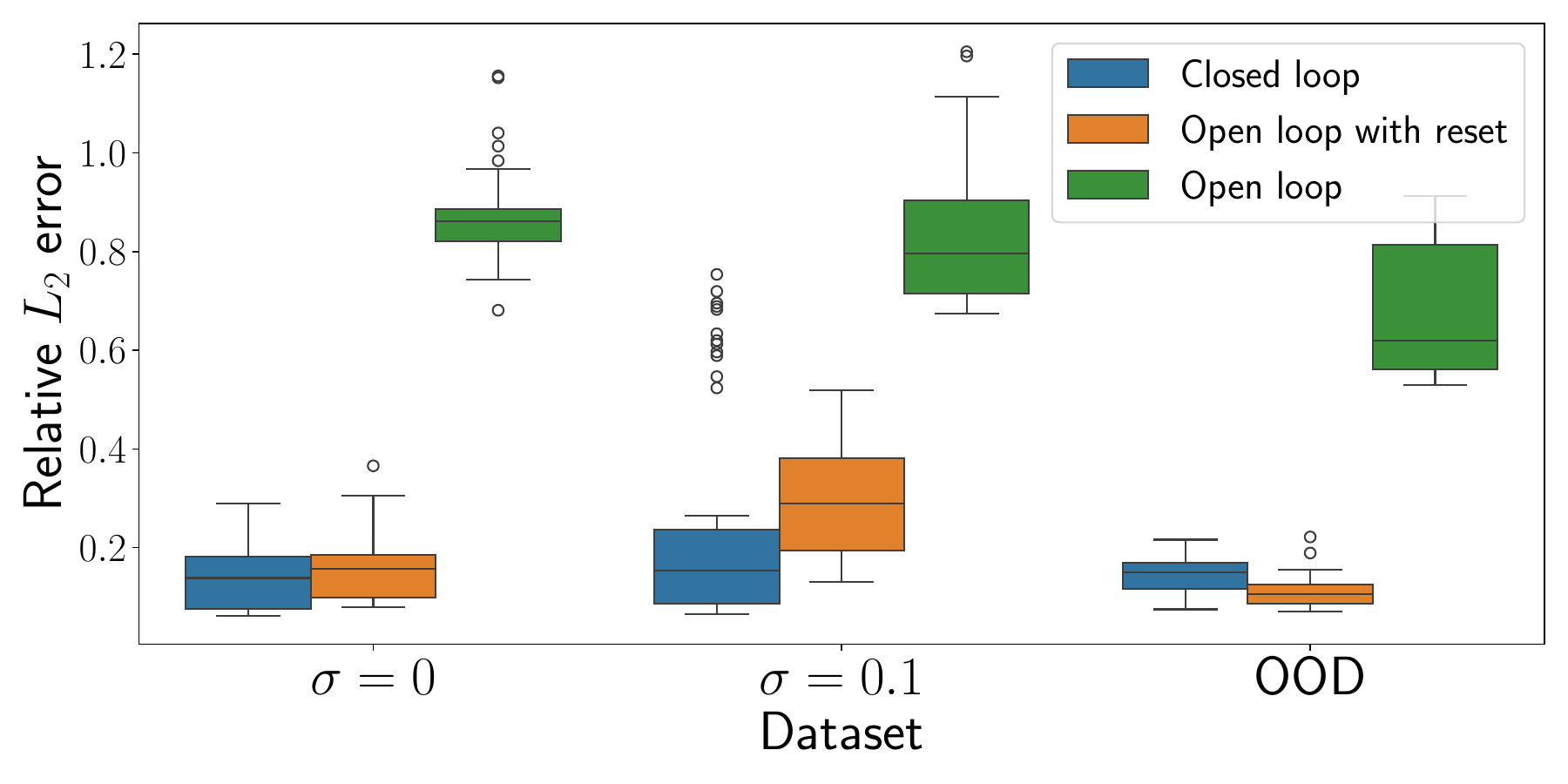}
    \caption{Robustness of the observers under noise and disturbance.}
    \label{fig:robustness}
    \vspace*{-0.3cm}
\end{figure}

The prediction accuracy of the observers across the test set in the noiseless, noisy, and OOD settings is shown in Fig. \ref{fig:robustness}. In the noiseless setting, the closed-loop observer $\opobscl$ generally performs better than the open-loop observer $\opobsol$ and the open-loop observer with reset $\opobsolr$. The poor performance of $\opobsol$ is explained by the instability caused by disturbance in the initial state. The improved performance of $\opobscl$ over $\opobsolr$ can be explained by the fact that $\opobscl$ benefits from the full historical dataset, as previous corrections influence all subsequent predictions through the continuous chain of autoregression. 

Moreover, $\opobscl$ is generally robust to noise, as the median error remains around the same level when the observer is applied to noisy data. This is not the case of $\opobsolr$, which sees a significant downgrade in performance. This can be explained by the high dependence on a few measurements in $\opobsolr$, rendering it vulnerable to noise. Instead, $\opobscl$ weighs the recent measurements against the prediction suggested by past measurements through the correction operator, resulting in a more accurate estimate. Furthermore, the accuracy of $\opobscl$ generally remains the same in the OOD dataset, suggesting that it is robust to distribution shifts. Overall, our simulations indicate that $\opobscl$ is both more accurate and more resilient to noise compared to $\opobsol$ and $\opobsolr$, as well as robust to distributional changes.

\begin{figure}[ht!]
    \centering
    \includegraphics[scale=0.28]{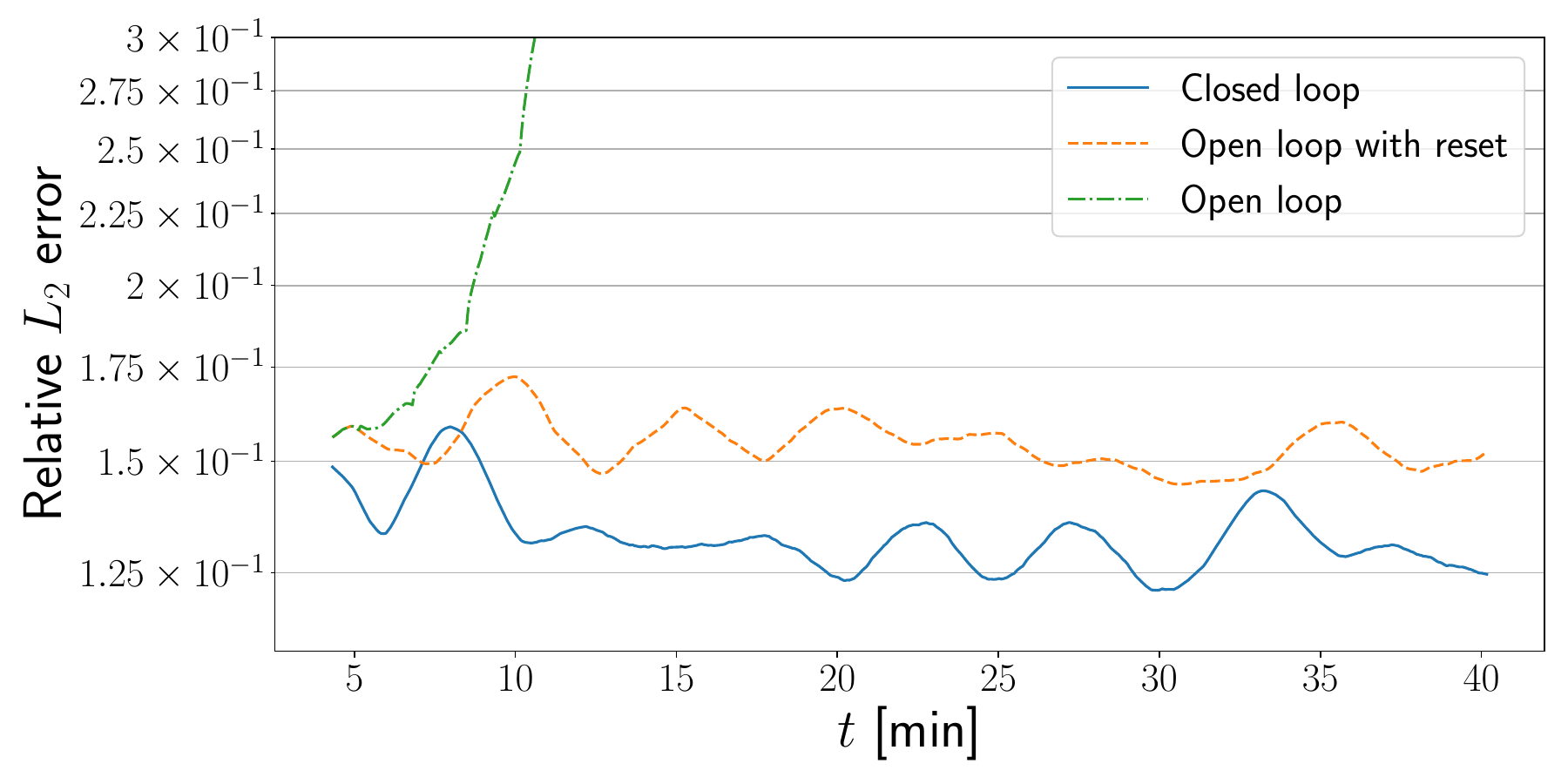}
    \caption{Evolution of prediction accuracy over the observation horizon.}
    \label{fig:l2_error_evolution}
    \vspace*{-0.5cm}
\end{figure}

Finally, the input-to-state stability for the observers is analyzed. In Fig. \ref{fig:l2_error_evolution}, the evolution of the prediction accuracy is compared across the observers as they are rolled out. The error grows unboundedly for the open-loop observer $\opobsol$, while it remains stable for the open-loop observer with reset $\opobsolr$, and gradually decreases for the closed-loop observer $\opobscl$. By construction, $\opobsol$ is completely dependent on an imperfect estimate of the initial state, which gives rise to the unbounded growth of the error. Meanwhile, $\opobsolr$ performs the same operation iteratively, so the error naturally remains stable. In contrast, $\opobscl$ cumulatively integrates data, resulting in an increasingly informed prediction manifested by a decreasing error. Our simulations thus indicate that $\opobscl$ exhibits the classical closed-loop property of ultimate-boundedness of the error.

\section{Conclusion}\label{sec:conclusion}
In this work, we introduced a data-driven, closed-loop observer for traffic density estimation using sparse stationary sensor measurements. Our method trains a Fourier neural operator to learn macroscopic traffic flow dynamics from high-fidelity data and integrates it with a correction operator to update predictions with real-time data. Numerical simulations with SUMO show that the closed-loop observer outperforms open-loop alternatives in terms of accuracy, robustness to noise, and stability over time, and indicate that the closed-loop observer is robust to distribution shifts. 

This study highlights the potential of this kind of learning-based, closed-loop observer and opens several interesting research directions. Future work includes leveraging real traffic data -- possibly with varying resolution -- to learn the dynamics, analyzing sensor placement for observability, and conducting formal analyses of convergence and stability.

\vspace*{-0.1cm}

\bibliographystyle{IEEEtran}
\bibliography{references}

\end{document}